\documentclass[12pt]{article}
\usepackage{amssymb, url, amsthm, amsmath}
\newtheorem{theorem}{\noindent Theorem}
\newtheorem{lemma}{\noindent Lemma}
\newtheorem{definition}{\noindent Definition}
\newtheorem{corollary}{\noindent Corollary}

\author{A.~Vershik\footnote{Supported by RFFI (grant no.~05-01-00899) and
CRDF (grant no.~RUM1-2622-ST-04).},
M.~Vsemirnov\footnote{Supported by the programs ``Scientific
Schools'' (grant no. NSh-8464-2006-1) and ``Modern problems of
theoretical mathematics''.}}
\date{11.04.07}
\title{The local stationary presentation of the alternating groups
and the normal form.}

\begin{document}
 \maketitle

\centerline{ABSTRACT}

We give the canonical normal form for the elements of the finite or
infinite alternating groups using local stationary presentation of
these groups.
 \section{The problem.}

The set of finite or infinite generators $x_1,\dots x_n$ ($n\in
{\mathbb N}$ or infinite) of a group $G$ is called the set {\it of
local generators of the depth $k$} if the following relations take
place:
  $$
  x_i\cdot x_j=x_j\cdot x_i,
  $$
when $|i-j|>k-1$ ($k=1$ means that $G$ is abelian). If the set of
relations between the elements $x_i \dots x_{i+k-1}$ do not depend
on $i$ (which means that the map sending $x_j \rightarrow
x_{j+i-1},\quad j=1, \dots, k$ is  an isomorphism between the
subgroup generated by $x_1,\dots x_k$ and the subgroup generated
by $x_i,\dots x_{i+k-1}$), we call the set of generators as {\it
stationary}. If this is true for all $i>i_0$ for some $i_0$ we
call it {\it eventually stationary}; see \cite{V1,V2,VO,V3}. We
will say in these cases that {\it the group has a local,
stationary (etc.) presentation}. The classical examples of the set
of local generators and local presentations for groups are Coxeter
generators and presentations of the Coxeter groups of infinite
series, standard generators of braid groups \cite{CM},
Curtis-Steinberg-Tits generators for classical matrix groups over
finite fields \cite[Section 4.1]{BGKLP}, etc.

The general problem, which appeared, is to describe the finite
groups which have the local stationary or eventually stationary
generators, to find such generators for given group if exists, to
find minimal possible depth $k$, etc.
 The problem could be considered in a wider context; for example for
 other parametrization of the generators that intervals of integers.
 Namely, instead of positive integers we can
consider the numeration of the generators by integers $\mathbb Z$,
or lattice ${\mathbb Z}^d$, or semilattice $\mathbb{N}^d$ or its
intervals, etc. The most interesting here are infinite locally
finite groups with infinite stationary set of local generators
numerated by the elements of a countable group or semigroup; like
f.e. infinite symmetric groups $S_\mathbb{N}$.

{\it In this paper we give the set of local stationary generators or
the depth 3 for the alternating group $A_n$} similar to the set of
of the depth 2 of classical Coxeter generators for symmetric groups.
The well-known classical presentations of the alternating groups
(see \cite {Car1,CM}) are very different and do not satisfy locality
conditions. In order to prove the main theorem {\it we give the
canonical normal form of the elements of the alternating group} with
respect to those generators.

\section{The Results}

 Consider the group generated by the generators $x_1,
\dots x_{n-2}$ subject to the following relations:
 \begin{eqnarray}
 \label{eq:relR}
  && x_i^3=1, \quad i=1, \dots, n-2, \\
 \label{eq:relS}
  && (x_i\cdot x_{i+1})^2=1 \quad i=1, \dots, n-3, \\
 \label{eq:relQ}
  && x_i\cdot x_j=x_j\cdot x_i,\quad i,j=1, \dots n-2, \quad
   |i-j|>2, \\
 \label{eq:relT}
  && x_i\cdot x_{i+1}^{-1}\cdot x_{i+2}= x_{i+2}\cdot x_i, \quad
  i=1 \dots n-4;
 \end{eqnarray}
   here  $n$ is either an integer greater than 1 or infinity.
We denote the free group with these relations by $S^+_n, n \in
\mathbb{N}$ (or $S^+_{\infty})$.

An equivalent form of the relation (\ref{eq:relT}) is  the
following:
 \begin{equation}
  \label{eq:relT'}
    x_{i+1}=[x_{i+2},x_i^{-1}] =
    x_{i+2} \cdot x_i^{-1} \cdot x_{i+2}^{-1} \cdot x_i.
 \end{equation}

\begin{theorem}
Let $n\ge 5$ be an integer. Then $S_n^+ \simeq A_n$. In another words the relations
{\rm(\ref{eq:relR})--(\ref{eq:relT})} define the stationary local presentation of the alternating
group:
$$A_n=S^+_n= \langle x_1, \dots x_{n-2}\,|\,(1)-(4) \rangle,\quad n=5,\dots; \quad
A_4=S^+_4=\langle x_1,x_2|(1),(2)\rangle.
$$
\end{theorem}
The cases $n=2,3$ is trivial.

 The next theorem shows that the set of relations can be reduced - the relations are not independent:

\begin{theorem}
Relations {\rm (\ref{eq:relT})}  for all $i\ge 2$ follow from all
relations {\rm (\ref{eq:relR})--(\ref{eq:relQ})} and relation {\rm
(\ref{eq:relT})} for $i=1$ {\rm(}i.e., $x_1\cdot x_2^{-1}\cdot
x_3=x_3\cdot x_1$ {\rm)}.
\end{theorem}

The key step in the proof of the Theorem 1 will be a normal form for
the elements of $S^+_{n}$ which is of independent interest.

We introduce the following notation.
\begin{definition}
For $m=1$, $2$,\dots, and $j=0,1,\dots,m+1$, we introduce the
following elements of $S_m^+$:
 \begin{eqnarray*}
 && y_{m,0}=x_{m}\cdot x_{m-1} \cdots x_2 \cdot x_1^2,\\
 && y_{m,k}=x_{m}\cdot x_{m-1} \cdots x_k, \quad k=1,\dots m, \\
 && y_{m,m+1}=\mathrm{id}.
 \end{eqnarray*}
In particular, $y_{1,0}=x_1^2$, $y_{1,1}=x_1$,
$y_{1,2}=\mathrm{id}$.
\end{definition}

\begin{theorem}[Normal form]
For each $n\ge 3$ and each element $X \in S^+_{n}$ there exist
integers $k_1$,\dots,$k_{n-2}$ such that $0\le k_j\le j+1$,
$j=1,\dots,n-2$ and the element $X$  has a representation  of the
form:
  \begin{equation}
  \label{eq:normalform}
   X=y_{1,k_1}\cdot y_{2,k_2} \cdot y_{n-2,k_{n-2}}.
  \end{equation}
In particular, $x_{n-2}$ appears in that form at most once, and the
generator $x_{n-k}$ appears at most $k-1$ times, $k=2, \dots n-1$.
\end{theorem}

The choice of the generators above for alternating group $A_n$ is
the following: $x_i=(i,i+1,i+2), i=1 \dots n-2$. It is not
difficult to prove that in the alternation group $A_n$ for all $n
\ne 6, n>2$ this is a unique solution (up to conjugacy) of the
system of relations above, and two solutions in the case
$A_6$.\footnote{For symmetric group $S_n$ the system of the
Coxeter relations also has unique solution up to conjugacy for all
$n \ne 6$; in $S_6$ there are two non conjugacy solutions of those
relations. The reason of that in both cases is that
$Out(S_6)=Aut(S_6)/Inn(S_6)={\mathbb{Z}}/2$; for $n\ne 6  - Aut
(S_n)=Inn(S_n)$.}
 Combining theorem 3 with theorem 1 we obtain the canonical normal
form for the elements of alternating groups with respect to those
generators.

\section{Proofs}

We start with the proof of the theorem 3 about normal form in the
group $S^+_{n}$.  Then we use it for the proof of theorem 1.

\begin{proof}[Proof of Theorem 3]
The proof of the theorem is similar to the the deduction  of the
canonical form for the elements of the symmetric group as Coxeter
group. It goes by induction on $n$ with the base $n=3,4$.

We will use the following transformation rules, which are
consequences of relations (\ref{eq:relR})--(\ref{eq:relT}):
  \begin{eqnarray}
   \label{eq:transf1}
   && x_{i+1}^2=x_{i+1}^{-1} = x_i\cdot x_{i+1}
   \cdot x_i, \quad
    i=1,2,\dots,n-3, \ \textrm{(relation~\ref{eq:relS})}, \\
   &&
   \label{eq:transf2}
   x_{i+2}\cdot x_i = x_i\cdot x_{i+1}^{-1}\cdot x_{i+2}
    \quad i=1,2,\dots,n-4, \ \textrm{(relation~\ref{eq:relT})} \\
   &&
   \label{eq:transf3}
    x_j \cdot x_i =  x_i\cdot x_j, \quad j-i\ge 3
    \ \textrm{(relation~\ref{eq:relQ})}.
  \end{eqnarray}

$\underline{n=3}$. In that case any element of $S_3^{+}$ can be
written in a unique way as was done above:
  $$
    x_1=y_{1,1},\ x_1^2=y_{1,0}, \ \mathrm{id}=y_{1,2}.
  $$
So, this is the group of order 3 isomorphic to $A_3$.

$\underline{n=4}$. Using (\ref{eq:relR}) we can write $X$ as a word
in $x_1$, $x_1^2$, $x_2$, $x_2^2$. We replace each occurrence of
$x_2^2$ by $x_1 \cdot x_2 \cdot x_1$ using transformation rule
(\ref{eq:transf1}) for $i=1$. Thus, we may assume that $x_2^2$ does
not appear in $X$. Next, any substring $x_2 \cdot x_1 \cdot x_2$ can
be replaced by $x_1^{2}=x_1^{-1}$  and any substring $x_2 \cdot
x_1^2 \cdot x_2$ can be transformed into
 $$
  x_2\cdot x_1 \cdot x_1 \cdot x_2 =
  x_1^{-1} \cdot x_2^{-1} \cdot  x_2^{-1} \cdot x_1^{-1} =
  x_1^2 \cdot x_2 \cdot x_1^2
 $$
(both transformations above are consequences of (\ref{eq:relS})
and (\ref{eq:relR})). Thus, we may assume that $x_2$ appears at
most once. In other words, any $X$ can be written in the form
 $$
  x_1^{a}, a=0,1,2;
  $$
  (this is the subgroup ${\mathbb Z}/3 \simeq A_3)$
  or
  $$ x_1^{a} x_2 x_1^{b},
  \textrm{ where } a,b \in \{0,1,2\},
 $$
which gives another 9 elements, so we have 12 elements of $S^+_4
\simeq A_4$.

\underline{Inductive step}. Assume the claim is true for $n\ge 4$,
we prove it for $n+1$.

Clearly, the subgroup of $S_{n+1}^+$ generated by the first $n-2$
letters $x_1$,\dots,$x_{n-2}$ is a quotient of $S_n^+$ (actually,
as we see later, they are isomorphic, but we do not use this fact
here; see Corollary~\ref{cor:2}). In particular, by the inductive
hypothesis any word in $x_1$,\dots,$x_{n-2}$ can be reduced to its
normal form (\ref{eq:normalform}).

First of all, each element  $X \in S^+_{n+1}$ can be represented
in the form:
  $$
  X=X_1\cdot x_{n-1}^{\alpha_1} \cdot X_2\cdot
  x_{n-1}^{\alpha_2}\dots x_{n-1}^{\alpha_m} \cdot X_{m+1},
  $$
where $m\ge 0$, $\alpha_j\in \{1,2\}, j=1, \dots, m $, and $X_j$,
$j=1,\dots,m+1$ are words in $x_1,\dots,x_{n-2}$. Because $n\ge
4$, we have the relation
 \begin{equation}
 \label{eq:xn-1^2}
 x_{n-1}^2 \stackrel{\textrm{by~(\ref{eq:relR})}}{=}
 x_{n-1}^{-1} \stackrel{\textrm{by~(\ref{eq:relS})}}{=}
 x_{n-2}\cdot x_{n-1} \cdot x_{n-2}.
 \end{equation}
Therefore, it is enough to consider the case where $\alpha_j=1$
for all $j=1,\dots,m$:
  \begin{equation}
  \label{eq:Xm}
  X=X'_1\cdot x_{n-1} \cdot X'_2\cdot
  x_{n-1}\dots x_{n-1} \cdot X'_{m+1}.
  \end{equation}
Assume that $m\ge 2$. We show that we can transform the right-hand
side of (\ref{eq:Xm}) into another word with fewer occurrences of
$x_{n-1}$. Consider the fragment of the above word between two
consecutive generators $x_{n-1}$:
 $$
  X= \dots x_{n-1} \cdot X'_j\cdot
  x_{n-1}\dots
 $$
By the inductive hypothesis $X'_j$ can be written in the normal
form (\ref{eq:normalform}). In particular, $x_{n-2}$ appears in
that normal form at most once. Using transformation rules
(\ref{eq:transf2})--(\ref{eq:transf3}) we can shift $x_{n-1}$ at
the left-hand side until we reach $x_{n-2}$ (if any) or the next
$x_{n-1}$. In the latter case we have $x_{n-1}^2$ and use
(\ref{eq:xn-1^2}) to diminish the number of $x_{n-1}$'s. In the
former case we obtain
 $$
  X= \dots x_{n-1} \cdot y_{n-2,k} \cdot x_{n-1}\dots
 $$
for some $k$, $0\le k \le n-1$. Now, using (\ref{eq:relQ}) we
shift $x_{n-1}$ at the right-hand side to the left until we reach
$x_{n-3}$ (if any) or $x_{n-2}$. There are small differences in
further analysis for $n=4$ and $n>4$. First assume that $n>4$. The
above transformations lead us to one of the two following
substrings:

  $x_{n-1}\cdot x_{n-2} \cdot x_{n-1}$, which is equal to
  $x_{n-2}^2$ by (\ref{eq:transf1}), or

  $x_{n-1}\cdot x_{n-2}\cdot x_{n-3} \cdot x_{n-1}$, which is
  equal to
\begin{eqnarray*}
 x_{n-1}\cdot x_{n-2}\cdot x_{n-3} \cdot x_{n-1}
 &\stackrel{\textrm{by~(\ref{eq:relS})}}{=}&
 x_{n-2}^{-1} \cdot x_{n-1}^{-1} \cdot x_{n-3} \cdot x_{n-1} \\
 &\stackrel{\textrm{by~(\ref{eq:relT})}}{=}&
 x_{n-2}^{-1} \cdot x_{n-3} \cdot x_{n-1}^{-1} \cdot x_{n-2} \cdot
  x_{n-1} \\
 &\stackrel{\textrm{by~(\ref{eq:relS})}}{=}&
 x_{n-2}^{-1} \cdot x_{n-3} \cdot x_{n-1}^{-2} \cdot x_{n-2}^{-1} \\
 &\stackrel{\textrm{by~(\ref{eq:relR})}}{=}&
 x_{n-2}^{2} \cdot x_{n-3} \cdot x_{n-1} \cdot x_{n-2}^{2}.
\end{eqnarray*}
In both cases we diminish the number of $x_{n-1}$'s.

If $n=4$, then there is one more case to be considered, namely,
the substring $x_3 \cdot x_2\cdot x_1^2 \cdot x_3$. We have
$$
x_3 \cdot x_2\cdot x_1^2 \cdot x_3
 \stackrel{\textrm{by~(\ref{eq:relT})}}{=}
 x_3 \cdot x_3 \cdot x_1^{-1} \cdot x_3^{-1} \cdot x_1
 \cdot x_1^2 \cdot x_3
 \stackrel{\textrm{by~(\ref{eq:relR})}}{=} x_3^2 \cdot x_1^{-1}
 \stackrel{\textrm{by~(\ref{eq:xn-1^2})}}{=}
 x_2 \cdot x_3 \cdot x_2 \cdot x_1^{-1},
$$
again diminishing the number of $x_3$'s.

Thus, it is enough to consider (\ref{eq:Xm}) with $m\le 1$.

If $m=0$, i.e., $x_{n-1}$ does not occur in $X$, we may apply the
inductive hypothesis and obtain the normal form
$$
   X=
   y_{1,k_1}\cdot y_{2,k_2} \cdots y_{n-2,k_{n-2}} \cdot \mathrm{id} =
    y_{1,k_1}\cdot y_{2,k_2} \cdots y_{n-2,k_{n-2}} \cdot y_{n-1,n}.
$$

If $m=1$, i.e., $X=X'_1 \cdot x_{n-1} \cdot X'_2$, we apply the
inductive hypothesis to $X'_2$ and write
 $$
   X=X'_1 \cdot x_{n-1} \cdot y_{1,j_1} \cdots y_{n-2,j_{n-2}}.
 $$
Using transformation rules (\ref{eq:transf2})--(\ref{eq:transf3})
we shift $x_{n-2}$ to the right until we reach $y_{n-2,j_{n-2}}$.
Therefore,
 $$
   X=X''_1 \cdot x_{n-1} \cdot y_{n-2,j_{n-2}}
   = X''_1 \cdot y_{n-1,j_{n-2}}.
 $$
Applying the inductive hypothesis to $X''_1$ we have
 $$
   X= y_{1,k_1}\cdot y_{2,k_2} \cdots y_{n-2,k_{n-2}}
   \cdot y_{n-1,j_{n-2}}.
 $$
thus completing the proof.
\end{proof}

\begin{lemma}
 \label{lem:lem1}
Consider the following elements of $A_n$: 
$x_i=(i,i+1,i+2), i=1, \dots n-2$. The relations {\rm
(\ref{eq:relR})--(\ref{eq:relT})} are true for these elements.

\end{lemma}
This fact is can be checked directly.

Thus, the group $A_n$ is a factor-group of the group $S^+_n$. To
prove the Theorem 1 it is enough to find the order of the group
$S^+_n$, namely, to prove that it contains at most $\frac{1}{2}n!$
elements.

\begin{proof}[Proof of Theorem~1]

 By Theorem 3 the order of $S_n^{+}$ is at most $\frac{1}{2} n!$.
Since $S_n^{+}$ projects onto $A_n$ by Lemma~\ref{lem:lem1}, we
conclude that $S_n^{+}\simeq A_n$.
\end{proof}

 In particular, this implies that for any
$X\in S_n^{+}\simeq A_n$ the normal form of the shape
(\ref{eq:normalform}) is unique.
\begin{proof}[Proof of Theorem 2] It follows by induction from the
calculations below. Suppose that relation (\ref{eq:relT}) (or its
equivalent form (\ref{eq:relT'})) is true for some $i$.  We prove it
for $i+1$. Namely,
\begin{eqnarray*}
     x_{i+3}\cdot x_{i+1}^{-1} \cdot x_{i+3}^{-1} \cdot
     x_{i+1} \cdot x_{i+2}^{-1}
      & \stackrel{\textrm{by~(\ref{eq:relR})}}{=} &
     x_{i+3} \cdot x_{i+1}^{-1}\cdot x_{i+3}^{-1}\cdot x_{i+1}
     \cdot x_{i+2} \cdot x_{i+2} \\
     & \stackrel{\textrm{by~(\ref{eq:relS})}}{=} &
     x_{i+3} \cdot x_{i+1}^{-1} \cdot x_{i+3}^{-1} \cdot
     x_{i+2}^{-1} \cdot x_{i+1}^{-1} \cdot x_{i+2} \\
    &\stackrel{\textrm{by~(\ref{eq:relS})}}{=}&
      x_{i+3} \cdot x_{i+1}^{-1} \cdot x_{i+2} \cdot x_{i+3}
     \cdot x_{i+1}^{-1} \cdot x_{i+2} \\
    &\stackrel{\textrm{by~(\ref{eq:relT'})}}{=}&
     x_{i+3} \cdot x_i^{-1} \cdot x_{i+2}\cdot x_i
     \cdot x_{i+3} \cdot x_i^{-1} \cdot x_{i+2} \cdot x_i \\
    &\stackrel{\textrm{by~(\ref{eq:relQ})}}{=}&
     x_{i+3} \cdot x_i^{-1} \cdot x_{i+2}
     \cdot x_{i+3} \cdot x_{i+2}\cdot x_i \\
    & \stackrel{\textrm{by~(\ref{eq:relQ})}}{=}&
      x_i^{-1} \cdot x_{i+3}\cdot x_{i+2}
     \cdot x_{i+3} \cdot x_{i+2} \cdot x_i \\
     & \stackrel{\textrm{by~(\ref{eq:relS})}}{=}&
      x_i^{-1} \cdot x_i =1.
\end{eqnarray*}
\end{proof}

We give two important corollaries of Theorems~1 and 3.
 Since the order of
$S_{n}^{+}=\frac{1}{2}n!$ we immediately have
\begin{corollary}
For each element $X\in S_n^{+}$ its normal form {\rm
(\ref{eq:normalform})} is unique.
\end{corollary}

\begin{corollary}
 \label{cor:2}
The subgroup of $S_{n+1}^{+}$ generated by the first $n-2$
generators $x_1$,\dots,$x_{n-2}$ is isomorphic to $S_n^+$.
\end{corollary}

\begin{proof}
Let $H$ be the subgroup of $S_{n+1}^{+}$ generated
$x_1$,\dots,$x_{n-2}$. Clearly, it is a factor of $S_n^+$. It
follows from the proof of Theorem~3, that $H$  has index $n+1$ in
$S_{n+1}^+$. By order considerations, $H$ must be isomorphic to
$S_n^+$ and we can identify these two groups.
\end{proof}

 Theorems 1 and 3 give the
way for the construction of the whole theory for alternating groups
independently from symmetric groups --- Bruhat order,
Gelfand-Tsetlin algebra and so on.
\section{The classical generators}
We conclude the paper by relating our presentation of the group
$A_n$ with the well-known one studied by Carmichael \cite{Car1}, see
modern explanation in \cite{CM}:
 \begin{equation}
  \label{eq:car}
  A_n \cong \langle v_1,\dots ,v_{n-2} : v_i^3=1, \ (v_i v_j)^2=1,
  \ i,j=1,\dots n-2, \ i\neq j \rangle.
 \end{equation}
Let us define $v_i$, $i=1,\dots,n-2$ by
 \begin{equation}
  \label{eq:defvi}
  v_i = \biggl(\prod_{k=i}^{n-2} x_k^{-1}\biggr)^{-1}
          \cdot \prod_{k=i+1}^{n-2} x_k^{-1}
      = x_{n-2} \cdots x_{i+1}\cdot x_i\cdot x_{i+1}^{-1}
      \cdots x_{n-2}^{-1}.
 \end{equation}
Using (\ref{eq:relR})--(\ref{eq:relT}) we can find that the
converse transformation is given by
\begin{equation}
  \label{eq:vtox}
  x_i = \biggl( \prod_{j=i+1}^{n-2} v_j \biggr)^{-1}
        \cdot \prod_{j=i}^{n-2} v_j =
        v_{n-2}^{-1}\cdots v_{i+1}^{-1} \cdot v_i \cdot v_{i+1}
        \cdots v_{n-2}.
\end{equation}
Thus $v_1$,\dots,$v_{n-2}$ also generate $S_n^+$. Moreover,  they
satisfy all the identities in presentation (\ref{eq:car}). This
gives an independent proof of Carmichael's result. Vice versa, one
can deduce our Theorem 1 from Carmichael's result. We leave details
for the reader.

\textbf{Remark.} In the recent paper \cite{R} \footnote{We are grateful to professor A.Postnikov
who had informed about paper \cite{R}} presentations of the analogues of the alternating groups was
given for all classical series of the Coxeter groups. In particular for the usual alternating group
the authors have used a well-known (see \cite{CM}) set of generators:
$$r_1=(2,3)(1,2)=(1,2,3),\ r_i=(1,2)(i+1,i+2) ,\ i>1;$$
and relations:
 \begin{eqnarray*}
 & r_1^3=r_i^2=1, \  i>1; \quad (r_1^{-1}r_2)^3=(r_i r_{i+1})^3=1,
  \ i=2,\dots n-2; &
 \\
 & (r_i r_j)^2=1,\ |i-j|>1. &
 \end{eqnarray*}
This set of the generators {\it is not local} because the last relation means that $r_i$ and $r_j$
commute only if the elements $r_i$ and $r_j$ are of order two, but $r_1^3=1\ne r_1^2$, so $r_1\cdot
r_i \ne r_i\cdot r_1$ for all $i>1$. It is interesting to find the sets of local generators for all
groups which was considered in  \cite{R}.


\begin{thebibliography}{6}
%
\bibitem{BGKLP} L.Babai, A.J. Goodman, W.M. Kantor, E.W. Luks and P.P.
P{\'a}lfy. Short presentations for finite groups. J. Algebra {\bf
194} (1997),
%
\bibitem{R} F.Brenti, V.Reiner, Y.Roihman. Alternating subgroups of
Coxeter groups. Nankai University. Tianjin. Conference FPSAC-2007.
Abstract.
%
\bibitem{Car1} R.D. Carmichael. Introduction to the theory of groups
of finite order. Dover Publications, Inc., New York, 1956.
%
\bibitem{CM} H.S.M. Coxeter, W.O.J. Moser. Generators and relations
for discrete groups. Springer-Verlag,
Berlin-Gottingen-Heidelberg, 1957.
%
\bibitem{V1} A.Vershik. Local stationary algebras. In ``Algebra and Analysis''. First
Siberian Winter School (Kemerovo, 1988). Amer. Math. Soc. Transl.
Ser. 2 {\bf 148}  (1991), 1--13.

%
\bibitem{V2} A.Vershik. Local algebras and a new version of Young's orthoganal form. In
``Topics in Algebra, part 2: Commutative rings and algebraic
groups'' (Warsaw 1988), Banach Cent. Publ. {\bf26}, Part 2,
(1990), 467--473.
%
\bibitem{VO}A.Okounkov, A.Vershik. A new approach to representation theory of symmetric
    groups. Selecta Math. {\bf 2} (1996), No.~4,  581--605.
%
\bibitem{V3}A.Vershik. Dynamic theory of growth in groups: entropy, boundaries,
examples. Russian Math. Surveys {\bf55} (2000), No.~4, 667--733.
\end{thebibliography}
\end{document}